# Modeling, Structure and Discretization of Mixed-dimensional Partial Differential Equations

J. M. Nordbotten and W. M. Boon
Department of Mathematics
University of Bergen
Norway

## Abstract

Mixed-dimensional partial differential equations arise in several physical applications, wherein parts of the domain have extreme aspect ratios. In this case, it is often appealing to model these features as lower-dimensional manifolds embedded into the full domain. Examples are fractured and composite materials, but also wells (in geological applications), plant roots, or arteries and veins.

In this manuscript, we survey the structure of mixed-dimensional PDEs in the context where the sub-manifolds are a single dimension lower than the full domain, including the important aspect of intersecting sub-manifolds, leading to a hierarchy of successively lower-dimensional sub-manifolds. We are particularly interested in partial differential equations arising from conservation laws. Our aim is to provide an introduction to such problems, including the mathematical modeling, differential geometry, and discretization.

## 1. Introduction

Partial differential equations (PDE) on manifolds are a standard approach to model on high-aspect geometries. This is familiar in the setting of idealized laboratory experiments, where 1D and 2D representations are used despite the fact that the physical world is 3D. Similarly, it is common to consider lower-dimensional models in applications ranging from geophysical applications. Some overview expositions for various engineering problems can be found in [1, 2, 3].

Throughout this paper we will consider the ambient domain to be 3D, and our concern is when models on 2D submanifolds are either coupled to the surrounding domain, and/or intersect on 1D and 0D submanifolds. Such models are common in porous media, where the submanifolds may represent either fractures (see e.g. [4]) or thin porous strata (see [1]), but also appear in materials [3]. In all these examples, elliptic differential equations representing physical conservation laws are applicable on all subdomains, and the domains of different dimensionality are coupled via discrete jump conditions. These systems form what we will consider as mixed-dimensional elliptic PDEs, and we will limit the exposition herein to this case.

In order to establish an understanding for the physical setting, we will in section two present a short derivation of the governing equations for fractured porous media, emphasizing the conservation structure and modeling assumptions. This derivation will lead to familiar models from literature (see e.g. [4, 5, 6, 7] and references therein).

We develop a unified treatment of mixed-dimensional differential operators on submanifolds of various dimensionality, using the setting of exterior calculus, and thus recast the physical problem in the sense of differential forms. We interpret the various subdomains as an imposed structure on the original domain, and provide a decomposition of differential forms onto the mixed-dimensional structure. By introducing a suitable inner product, we show that this mixed-dimensional space is a Hilbert space. On this decomposition we define a semi-discrete exterior derivative, which leads to a de Rham complex with the same co-homology structure as the original domain. A co-differential operator can be defined via the inner product, and it is possible to calculate an explicit expression for the co-differential operator. This allows us to establish a Helmholtz decomposition on the mixed-dimensional geometry. We also define the mixed-dimensional extensions of the familiar Sobolev spaces.

Having established the basic ingredients of a mixed-dimensional calculus, we are in a position to discuss elliptic minimization problems. Indeed, the mixed dimensional minimization problems are well-posed with unique solutions based on standard arguments, and we also state the corresponding Euler equations (variational equations). With further regularity assumptions, we also give the strong form of the minimization problems, corresponding to conservation laws and constitutive laws for mixed-dimensional problems.

## 2. Fractured porous media as a mixed-dimensional PDE

This section gives the physical motivation for mixed-dimensional PDE. As the section is meant to be motivational, we will omit technical details whenever convenient. We will return to these details in the following sections.

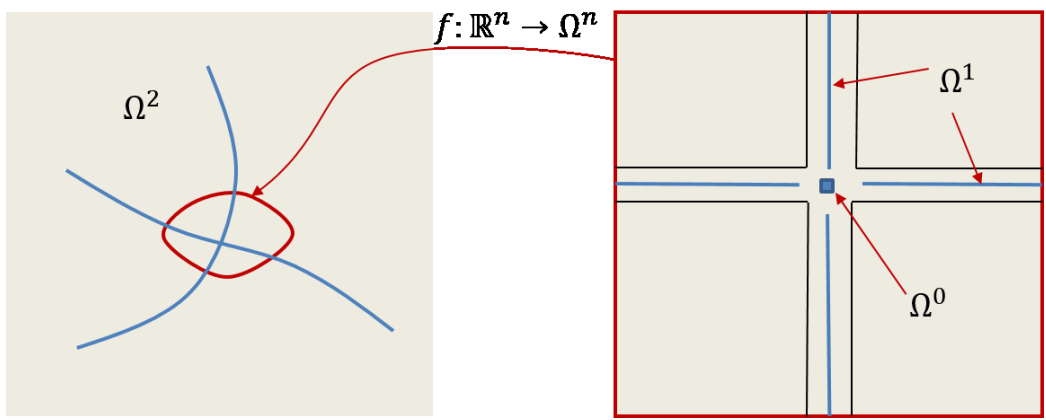

**Figure 1**: Example geometry of two intersecting fractures in 2D, and the logical representation of the intersection after mapping to a local coordinate system.

We consider the setting of a domain $D \in \mathbb{R}^n$. In sections 3 and onwards we will consider arbitrary $n$, however in this section we will for simplicity of exposition consider only $n = 3$. We consider a fractured media, where we are given explicit knowledge of the fractures, thus we consider the domains $\Omega_i^d$ as given, where $d$ represents the dimensionality of the domain and $i$ is an index. In particular, intact material lies in domains of $d = 3$, while $d = 2$ represents fracture segments, and $d = 1$ represents intersections, see Figure 1. For each domain $\Omega_i^d$ we assign an orientation based on $n-d$ outer normal vectors $\boldsymbol{\nu}_{ij}$.

In order to specify the geometry completely, we consider the index sets $\hat{S}_i$ and $\check{S}_i$ as the $d+1$ dimensional and $d-1$ dimensional neighbors of a domain $i$. Thus for $d=2$, the set $\hat{S}_i$ contains the domain(s) $\Omega_l^3$ which are on the positive (and negative) side of $\Omega_i^2$. On the other hand, the set $\check{S}_i$ contains the lines that form (parts of) the boundary of $\Omega_i$. Additionally, the set of all lower-dimensional neighbors is defined as $\check{\mathfrak{S}}_i = \left[\check{S}_i, S_{\check{S}_i}, \dots\right]$ We will use a summation convention over omitted indexes, thus $\Omega^d$ represents all subdomains of dimension $d$, while $\Omega$ is the full mixed-dimensional stratification.

For steady-state flows in porous media, the fluid satisfies a conservation law, which for intact rock and an $n$-dimensional fluid flux vector $\boldsymbol{u}$ takes the form

$$\nabla \cdot \boldsymbol{u} = f \qquad \text{on} \qquad D \qquad (2.1)$$

We wish to express this conservation law with respect to our geometric structure. To this end, let us first define the mixed-dimensional flux $\mathbf{u}$, which is simply a $d$-dimensional vector field on each $\Omega_i^d$. We write $\mathbf{u} = \left[\boldsymbol{u}_i^d\right]$ when we want to talk about specific components of $\mathbf{u}$. We similarly define other mixed-dimensional variables, such as the source-term $\mathsf{f}$.

Now clearly, for $d=n$, we recover equation (2.1). Now consider $d=n-1$, and a fracture of variable Lipschitz-continuous aperture as indicated in figure 2.

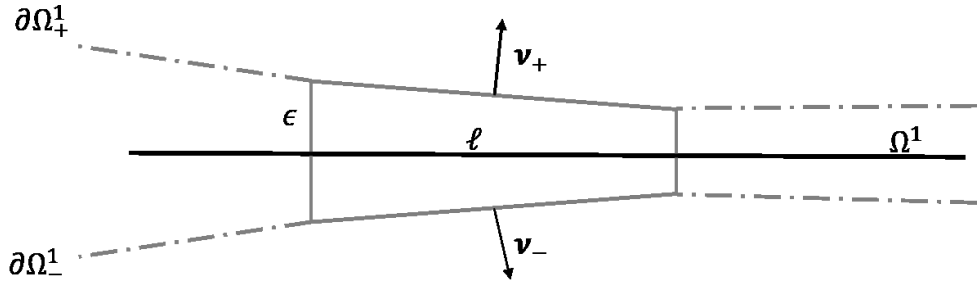

**Figure 2**: Example of local geometry for derivation of mixed-dimensional conservation law.

Here the dashed lines indicate a fracture boundary, the solid black line is the lower-dimensional representation, and the solid gray line indicates the region of integration, $\omega$, of length $\ell$ and width $\epsilon(x)$. Evaluating the conservation law over $\omega$ leads to

$$\int_\omega \nabla \cdot \boldsymbol{u} \; da = \int_{\partial\omega} \boldsymbol{u} \cdot \boldsymbol{\nu} \; ds = \int_\omega \phi$$

where $\boldsymbol{\nu}$ are the external normal vectors. Since our integration area is in the limiting case of $\ell \to 0$ a quadrilateral, we split the last integral into parts where $\boldsymbol{\nu}$ is constant,

$$\int_{\partial\omega} \boldsymbol{u} \cdot \boldsymbol{\nu} \; ds = \boldsymbol{\nu}_+ \cdot \int_{\partial\omega_+} \boldsymbol{u}_+^3 \; ds + \boldsymbol{\nu}_- \cdot \int_{\partial\omega_-} \boldsymbol{u}_-^3 \; ds + \int_{\partial\omega_r} \boldsymbol{\tau} \cdot \boldsymbol{u} \; ds - \int_{\partial\omega_l} \boldsymbol{\tau} \cdot \boldsymbol{u} \; ds$$

Here we denote the various parts of the boundary as + and − (for up and down) and l̲eft and r̲ight. The notation $\boldsymbol{\tau}$ is the tangential vector to $\Omega^1$. Clearly, letting the length $\ell$ be infinitesimal, the last two terms satisfy

$$\lim_{\ell \to 0} \frac{\int_{\partial \omega_r} \boldsymbol{\tau} \cdot \boldsymbol{u}\ ds - \int_{\partial \omega_l} \boldsymbol{\tau} \cdot \boldsymbol{u}\ ds}{\ell} = \nabla_{\Omega^2} \cdot \int_{\omega_-}^{\omega_+} \boldsymbol{\tau} \cdot \boldsymbol{u}\ ds = \nabla_{\Omega^2} \cdot (\epsilon \boldsymbol{u}^2)$$

where $\nabla_{\Omega^2} \cdot$ is the in-plane divergence and

$$\boldsymbol{u}^2 \equiv \frac{1}{\epsilon} \int_{\omega_-}^{\omega_+} \boldsymbol{\tau} \cdot \boldsymbol{u}\ ds \tag{2.3}$$

Considering similarly the limits of $\ell \to 0$ for the two first terms, we obtain e.g.

$$\lim_{\ell \to 0} \boldsymbol{v}_+ \cdot \ell^{-1} \int_{\partial \omega_+} \boldsymbol{u}_+^3\ ds = \left( 1 + \left| \frac{d}{dx} \nabla_{\Omega^2}(\partial \omega_+) \right|^2 \right)^{1/2} \boldsymbol{v}_+ \cdot \boldsymbol{u}_+^3$$

Combining the above, we thus have

$$\lim_{\ell \to 0} \ell^{-1} \int_\omega \nabla \cdot \boldsymbol{u}\ da = \lambda_+^3 + \lambda_-^3 + \nabla_{\Omega^2} \cdot (\epsilon \boldsymbol{u}^2) = [\![\lambda^3]\!] + \nabla_{\Omega^{2w}} \cdot (\epsilon \boldsymbol{u}^2)$$

$$\tag{2.4}$$

where $\lambda_\pm^3$ is defined as

$$\lambda_\pm^3 = \left( 1 + \left| \frac{d}{dx} \nabla_{\Omega^2}(\partial \omega_+) \right|^2 \right)^{1/2} \boldsymbol{v}_\pm^3 \cdot \boldsymbol{u}_\pm^3 \tag{2.5}$$

and

$$[\![\lambda^3]\!] = -\sum_{l \in \hat{S}_i} \lambda_l^3 \tag{2.5b}$$

Note that we have made no approximations in obtaining equation (2.4) – the left-hand side is an exact expression of conservation. The model approximations appear later when deriving suitable constitutive laws for i.e. $\boldsymbol{u}^2$, etc. In practice, since the fractures have a high aspect ratio by definition, the pre-factor in equation (2.5) is approximately identity, and it is common to simply use the approximation

$$\lambda_\pm^3 = \boldsymbol{v}_\pm^3 \cdot \boldsymbol{u}_\pm^3 \tag{2.5c}$$

The derivation above (including the definition in equation (2.4)), generalizes in the same way to intersection lines and intersection points, thus we find that for all $d < n$ it holds that

$$[\![\hat{\epsilon}\hat{\lambda}]\!] + \nabla_{\Omega^d} \cdot (\epsilon^d \boldsymbol{u}^d) = \phi^d \tag{2.6}$$

Here the hat again denotes the next higher-dimensional domains, so that $\hat{\lambda} = \lambda^{d+1}$. By defining $\lambda^{n+1} = 0$, it is clear that equation (2.6) also holds for $d = n$, and thus it represents the mixed-dimensional conservation law for all $\Omega_i^d$. In this more general setting, $\epsilon$ denotes the cross-sectional width (2D), area (1D) and volume (0D) for successively lower-dimensional intersections.

For porous materials, the conservation law (2.1) is typically closed by introducing Darcy's law as a modeling assumption, stated in terms of a potential $p$ on the domain $D$ as

$$\boldsymbol{u} = -k\nabla p \tag{2.7}$$

The coefficient $k$ is in general a tensor. Unlike for the conservation law, it is not possible to derive an exact expression for the mixed-dimensional constitutive law, but by making some (reasonable) assumptions on the structure of the solution, it is usually accepted that Darcy's law is inherited for each subdomain (see extended discussion in [1], but also [8]), i.e.

$$\boldsymbol{u}^d = -k^d \nabla_{\Omega^1} p^d \tag{2.8}$$

To close the model, it is also necessary to specify an additional constraint, where the two most common choices are that either the potential is continuous (see discussion in [9])

$$\hat{p}_+ = \hat{p}_- \tag{2.9}$$

or, more generally, that the pressure is discontinuous but related to the normal flux above

$$\hat{\lambda}_{\pm} = -2\hat{k}_{\pm}\frac{p^d - \hat{p}_{\pm}}{(\epsilon^d)^{\frac{1}{n-d}}} \tag{2.10}$$

The model equations (2.6), (2.8) together with (2.9) and (2.10) are typical of those used in practical applications. However, to the authors' knowledge, our work is first time they are explicitly treated as mixed-dimensional PDE (see also [10, 11]).

## 3. Exterior calculus for mixed-dimensional geometries

We retain the same geometry as in the previous section, but continue the exposition in the language of exterior calculus (for introductions, see [12, 13, 14]). First, we note that the components of the mixed-dimensional flux presented above all correspond to $d-1$ forms, $\boldsymbol{u}_i^d \in \Lambda^{d-1}(\Omega_i^d)$, while the components of pressure all correspond to $d$-forms, $p_i^d \in \Lambda^d(\Omega_i^d)$. This motivates us to define the following mixed-dimensional $k$-form

$$\mathfrak{L}^k(\Omega) = \prod_{i,d} \Lambda^{k-(n-d)}(\Omega_i^d) \tag{3.1}$$

From here on, it is always assumed that $\mathfrak{L}^k$ is defined over $\Omega$, and the argument is suppressed.

Moreover, we note that equation (2.5c) is (up to a sign) the trace with respect to the inclusion map of the submanifold, thus for a mixed-dimensional variable $\mathfrak{a} \in \mathfrak{L}^k$ the jump operator is naturally written as

$$(\mathbb{d}\mathfrak{a})_i^d = (-1)^{d+k} \sum_{j \in \hat{S}_i} \varepsilon(\Omega_i^d, \partial_i \Omega_j^{d+1}) \text{Tr}_{\Omega_i^d} a_j^{d+1} \tag{3.2}$$

Here we have exchanged the bracket notation of equation (2.5b), which is common in applications, with a simpler notation, $\mathbb{d}$, which more clearly emphasizes that this is a (discrete) differential operator, in the normal direction(s) with respect to the submanifold. We use the notation $\varepsilon(\Omega_i^d, \partial_i \Omega_j^{d+1})$ to indicate the relative orientation (positive or negative) of the arguments.

We obtain a mixed-dimensional exterior derivative, which we denote $\mathfrak{d}$, by combining the jump operator with the exterior derivative on the manifold,

$$(\mathfrak{d}\mathfrak{a})_i^d = da_i^d + (\mathbb{d}\mathfrak{a})_i^d \tag{3.3}$$

This expression is meaningful, since both $da_i^d$, $(\mathbb{d}a)_i^d \in \Lambda^{k-(n-d)+1}(\Omega_i^d)$, and thus clearly $\mathfrak{d}a \in \mathfrak{L}^{k+1}$. A straight-forward calculation shows that $d(\mathbb{d}a)_i^d = -(\mathbb{d}da)_i^d$, thus for all $a$

$$\mathfrak{d}\mathfrak{d}a = 0 \qquad (3.4)$$

and it can furthermore be shown that if $a = 0$, and if $D$ is contractible, then there exists $b \in \mathfrak{L}^{k-1}$ such that $a = \mathfrak{d}b$. Thus the mixed-dimensional exterior derivative forms a de Rham complex,

$$0 \to \mathbb{R} \overset{\subset}{\to} \mathfrak{L}^0 \overset{\mathfrak{d}}{\to} \mathfrak{L}^1 \overset{\mathfrak{d}}{\to} ... \overset{\mathfrak{d}}{\to} \mathfrak{L}^n \to 0 \qquad (3.5)$$

which is exact (for the proof of this, and later assertions, please confer [11]).

The natural inner product for the mixed-dimensional geometry must take into account the traces on boundaries, and thus takes the form for $a, b \in \mathfrak{L}^k$

$$(a, b) = \sum_{i,d}\big(a_i^d, b_i^d\big) + \sum_{l \in \tilde{\mathfrak{S}}_i^d}\big(\mathrm{Tr}_{\Omega_l} a_i^d, \mathrm{Tr}_{\Omega_l} b_i^d\big) \qquad (3.6)$$

It is easy to verify that this is indeed an inner product, and thus forms the norm on $\mathfrak{L}^k$

$$\|a\| = (a, a)^{1/2} \qquad (3.7)$$

The codifferential $\mathfrak{d}^*: \mathfrak{L}^k \to \mathfrak{L}^{k-1}$ is defined as the dual of the exterior derivative with respect to the inner product, such that for $a \in \mathfrak{L}^k$

$$(\mathfrak{d}^*a, b) = (a, \mathfrak{d}b) + (\mathrm{Tr}\, b, \mathrm{Tr}^*a)_{\partial D} \qquad \text{for all } b \in \mathfrak{L}^{k-1} \qquad (3.8)$$

It follows from the properties of inner product spaces that the codifferential also forms an exact de Rham sequence. Thus, when $D$ is contractible, we have the following Helmholtz decomposition: For all $a \in \mathfrak{L}^k$, there exist $a_{\mathfrak{d}} \in \mathfrak{L}^{k-1}$ and $a_{\mathfrak{d}^*} \in \mathfrak{L}^{k+1}$ such that

$$a = \mathfrak{d}a_{\mathfrak{d}} + \mathfrak{d}^*a_{\mathfrak{d}^*} \qquad (3.9)$$

In view of the uncertainty in the modeling community of the correct constitutive laws for mixed-dimensional problems (as per the discussion of equation (2.9) and (2.10)), it is of great practical utility to be able to explicitly calculate the co-differential, since this will have the structure of the constitutive law. Utilizing equations (3.6) and (3.8), we obtain

$$(\mathfrak{d}^*b)_i^d = d^*b_i^d \qquad \qquad \text{on } \Omega_i^d \qquad (3.10)$$

and

$$\mathrm{Tr}_{\partial_j \Omega_i^d}(\mathfrak{d}^*b)_i^d = \mathrm{d}^*\mathrm{Tr}_{\partial_j \Omega_i^d} b_i^d + \left((-1)^{d+k}\varepsilon\big(\Omega_j^{d-1}, \partial_j \Omega_i^d\big) b_j^{d-1} - \mathrm{Tr}^*_{\partial_j \Omega_i^d} b_i^d\right) \quad \text{on } \partial \Omega_i^d \qquad (3.11)$$

We close this section by noting that the differential operators provide the basis for extending Hilbert spaces to the mixed-dimensional setting. In particular, we are interested in the generalized $L_2$ space

$$L_2 \mathfrak{L}^k : \{a \in \mathfrak{L}^k \mid \|a\| < \infty\} \qquad (3.12)$$

as well as the first order differential spaces

$$H\mathfrak{L}^k : \{a \in L_2 \mathfrak{L}^k \mid \|\mathfrak{d}a\| < \infty\} \quad \text{and} \quad H^*\mathfrak{L}^k : \{a \in L_2 \mathfrak{L}^k \mid \|\mathfrak{d}^*a\| < \infty\} \qquad (3.13)$$

We denote the norms of $H\mathfrak{L}^k$ and $H^*\mathfrak{L}^k$ by

$$\|\mathfrak{a}\|_H = \|\mathfrak{a}\| + \|\mathfrak{d}\mathfrak{a}\| \qquad \text{and} \qquad \|\mathfrak{a}\|_{H^*} = \|\mathfrak{a}\| + \|\mathfrak{d}^*\mathfrak{a}\| \tag{3.14}$$

Before closing, we introduce the convention that a circle above the function space denotes homogeneous boundary conditions, i.e. $\mathring{H}\mathfrak{L}^k: \{\mathfrak{a} \in H\mathfrak{L}^k \,|\, \mathrm{Tr}_{\partial D}\,\mathfrak{a} = 0\}$ and $\mathring{H}^*\mathfrak{L}^k: \{\mathfrak{a} \in H^*\mathfrak{L}^k \,|\, \mathrm{Tr}^*_{\partial D}\,\mathfrak{a} = 0\}$.

Then, the Poincaré inequality holds for contractible domains in the mixed-dimensional setting for either $\mathfrak{a} \in \mathring{H}\mathfrak{L}^k \cap H^*\mathfrak{L}^k$ or $\mathfrak{a} \in H\mathfrak{L}^k \cap \mathring{H}^*\mathfrak{L}^k$:

$$\|\mathfrak{a}\| \leq C_\Omega(\|\mathfrak{d}\mathfrak{a}\| + \|\mathfrak{d}^*\mathfrak{a}\|) \tag{3.15}$$

## 4. Mixed-dimensional elliptic PDEs

Based on the extension of the exterior derivative and its dual to the mixed-dimensional setting, we are now prepared to define the generalization of elliptic PDEs. We start by considering the minimization problem equivalent to the Hodge Laplacian for $\mathfrak{a} \in \mathfrak{L}^k$

$$\mathfrak{a} = \arg \inf_{\mathfrak{a} \in \mathring{H}\mathfrak{L}^k \cap H^*\mathfrak{L}^k} J_{\mathfrak{K}}(\mathfrak{a}') \tag{4.1}$$

where we define the functional by

$$J_{\mathfrak{K}}(\mathfrak{a}') = \frac{1}{2}(\mathfrak{K}\mathfrak{d}^*\mathfrak{a}', \mathfrak{d}^*\mathfrak{a}') + \frac{1}{2}(\mathfrak{K}^*\mathfrak{d}\mathfrak{a}', \mathfrak{d}\mathfrak{a}') - (\mathfrak{f}, \mathfrak{a}') \tag{4.2}$$

For equation (4.1) to be well-posed and have a unique solution, we need $(\mathfrak{K}\mathfrak{d}^*\mathfrak{a}', \mathfrak{d}^*\mathfrak{a}') + (\mathfrak{K}^*\mathfrak{d}\mathfrak{a}', \mathfrak{d}\mathfrak{a}')$ to be continuous and coercive, i.e. we need to impose constraints on $\mathfrak{K}$ and $\mathfrak{K}^*$. Indeed, by reverting to the definition of the inner product, we have that

$$(\mathfrak{K}\mathfrak{d}^*\mathfrak{a}', \mathfrak{d}^*\mathfrak{a}') + (\mathfrak{K}^*\mathfrak{d}\mathfrak{a}', \mathfrak{d}\mathfrak{a}') \geq \min(\alpha_{\mathfrak{K}}, \alpha_{\mathfrak{K}^*})\,(1 + C_\Omega)^2\|\mathfrak{a}'\|^2$$

Where the ellipticity constant $\alpha_{\mathfrak{K}}$ is the minimum eigenvalue of $\mathfrak{K}$, and similarly for $\alpha_{\mathfrak{K}^*}$.

The minimum of equation (4.1) must satisfy the Euler-Lagrange equations, thus $\mathfrak{a} \in \mathring{H}\mathfrak{L}^k \cap H^*\mathfrak{L}^k$ satisfies

$$(\mathfrak{K}\mathfrak{d}^*\mathfrak{a}, \mathfrak{d}^*\mathfrak{a}') + (\mathfrak{K}^*\mathfrak{d}\mathfrak{a}, \mathfrak{d}\mathfrak{a}') = (\mathfrak{f}, \mathfrak{a}') \qquad \text{for all } \mathfrak{a}' \in \mathring{H}\mathfrak{L}^k \cap H^*\mathfrak{L}^k \tag{4.3}$$

From the perspective of applications, and mirroring the distinctions between conservation laws and constitutive laws discussed in Section 2, we will be interested in the mixed formulation of equation (4.3) obtained by introducing the variable $\mathfrak{b} = \mathfrak{K}\mathfrak{d}^*\mathfrak{a}$. Then we may either consider a constrained minimization problem derived from equation (4.1), or for the sake of brevity, proceed directly to the Euler-Lagrange formulation: Find $(\mathfrak{a}, \mathfrak{b}) \in H\mathfrak{L}^k \times H\mathfrak{L}^{k-1}$ which satisfy

$$(\mathfrak{K}^{-1}\mathfrak{b}, \mathfrak{b}') - (\mathfrak{a}, \mathfrak{d}\mathfrak{b}') = 0 \qquad \text{for all } \mathfrak{b}' \in H\mathfrak{L}^{k-1} \tag{4.4}$$

$$(\mathfrak{d}\mathfrak{b}, \mathfrak{a}') + (\mathfrak{K}^*\mathfrak{d}\mathfrak{a}, \mathfrak{d}\mathfrak{a}') = (\mathfrak{f}, \mathfrak{a}') \qquad \text{for all } \mathfrak{a}' \in H\mathfrak{L}^k \tag{4.5}$$

The saddle-point formulation is well-posed subject to Babuška-Aziz inf-sup condition. Due to the presence of a Helmholtz decomposition, this follows by standard arguments. From equations (4.4) and (4.5) we deduce the strong form of the Hodge Laplacian on mixed form, corresponding to the equations

$$\mathfrak{b} = \mathfrak{K}\mathfrak{d}^*\mathfrak{a} \qquad \text{and} \qquad \mathfrak{d}\mathfrak{b} + \mathfrak{d}^*(\mathfrak{K}^*\mathfrak{d}\mathfrak{a}) = \mathfrak{f} \tag{4.6}$$

Of the various formulations, equations (4.4) and (4.5) are particularly appealing as they do not require the coderivative.

An important remark is that the relative simplicity of the well-posedness analysis for the mixed-dimensional equations relies on the definition of the function spaces and norms. In particular, due to the definition of $H\mathfrak{L}^k$ via the mixed-dimensional differential $\mathfrak{d}$, the norm on the function space is inherently also mixed-dimensional, and cannot simply be decomposed into, say norms on the function spaces $H\Lambda^{k-(n-d)}(\Omega_i^d)$. For this reason, analysis in terms of "local norms" becomes significantly more involved [15, 16, 17].

## 5. Finite-dimensional spaces

In order to exploit the mixed-dimensional formulations from the preceding section, and in particular equations (4.4-4.5) we wish to consider finite-dimensional subspaces of $H\mathfrak{L}^k$. These spaces should be constructed to inherit the de Rham structure of equation (3.5), and with bounded projection operators. A natural approach is to consider the polynomial finite element spaces as a starting point [13].

From the finite element exterior calculus (FEEC - [13]), we know that on the highest-dimensional domains $\Omega_i^d$, we may choose any of the finite element de Rham sequences, and in particular, we may consider the standard spaces from applications for a simplicial tessellation $\mathcal{T}_i^n = \mathcal{T}(\Omega_i^n)$

$$\mathcal{P}_r\Lambda^k(\mathcal{T}_i^n) \qquad \text{and} \qquad \mathcal{P}_r^-\Lambda^k(\mathcal{T}_i^n) \tag{5.1}$$

These correspond to the full and reduced polynomial spaces of order $r$, respectively. In order to build a finite element de Rham sequence, we recall that (while still commuting with bounded projection operators) the full polynomial spaces satisfy

$$\mathcal{P}_r\Lambda^k(\mathcal{T}^n) \xrightarrow{d} \mathcal{P}_{r-1}\Lambda^{k+1}(\mathcal{T}^n) \qquad \text{and} \qquad \mathcal{P}_r\Lambda^k(\mathcal{T}^n) \xrightarrow{d} \mathcal{P}_r^-\Lambda^{k+1}(\mathcal{T}^n) \tag{5.2}$$

while the reduced spaces satisfy

$$\mathcal{P}_r^-\Lambda^k(\mathcal{T}^n) \xrightarrow{d} \mathcal{P}_r^-\Lambda^{k+1}(\mathcal{T}^n) \qquad \text{and} \qquad \mathcal{P}_r^-\Lambda^k(\mathcal{T}^n) \xrightarrow{d} \mathcal{P}_{r-1}\Lambda^{k+1}(\mathcal{T}^n) \tag{5.3}$$

Thus, any of these combination of spaces are acceptable for $\Omega_i^n$, and consider therefore the choice as given, and denoted by $\Lambda_h^{k,n}$ and $\Lambda_h^{k+1,n}$.

For $d < n$, we must consider not only the continuous differential operator $d$, but also the discrete jump operator $\mathbb{d}$. It is therefore clear that for i.e. $d = n - 1$, we must consider the traces of the spaces chosen for $\Omega^n$. In particular, we require for all pairs of dimensions $0 \leq e < d \leq n$,

$$\mathrm{Tr}_{\Omega_i^e}\,\Lambda_h^k(\mathcal{T}^d) \subseteq \Lambda_h^{k+(n-e)}(\mathcal{T}^e) \tag{5.4}$$

The traces of the standard finite element spaces can be summarized as follows for $e < d$ [13]:

$$\text{Tr}_{\Omega^e}\, \mathcal{P}_r \Lambda_h^k(\mathcal{T}^d) = \mathcal{P}_{r-e+k}^- \Lambda_h^{k+(n-e)}(\mathcal{T}^e) \qquad \text{and} \qquad \text{Tr}_{\Omega^e}\, \mathcal{P}_r^- \Lambda_h^k(\mathcal{T}^d) = \mathcal{P}_{r-e+k-1} \Lambda_h^{k+(n-e)}(\mathcal{T}^e)$$

(5.5)

We now define the polynomial subspaces $\mathcal{P}_{\mathfrak{r}}^{\mathfrak{m}} \mathfrak{L}^k \in H\mathfrak{L}^k$ as

$$\left(\mathcal{P}_{\mathfrak{r}}^{\mathfrak{m}} \mathfrak{L}^k\right)_i^d = \mathcal{P}_{r_i^d}^{p_i^d} \Lambda^{k-(n-d)}(\mathcal{T}_i^d)$$

(5.6)

where the multi-indexes $\mathfrak{r}$ and $\mathfrak{m}$ have values $r_i^d \in \mathbb{P}$ and $m_i^d \in [\,,-]$, respectively. When the multi-indexes are chosen to satisfy both (5.2-5.3) as well as (5.4), we obtain the discrete de Rham complex

$$0 \to \mathbb{R} \hookrightarrow \mathcal{P}_{\mathfrak{r}}^{\mathfrak{m}} \mathfrak{L}^0 \xrightarrow{\mathfrak{d}} \mathcal{P}_{\mathfrak{r}}^{\mathfrak{m}} \mathfrak{L}^1 \xrightarrow{\mathfrak{d}} \ldots \xrightarrow{\mathfrak{d}} \mathcal{P}_{\mathfrak{r}}^{\mathfrak{m}} \mathfrak{L}^n \to 0$$

(5.7)

Due to the existence of stable projections for all finite element spaces in $\mathcal{P}_{\mathfrak{r}}^{\mathfrak{m}} \mathfrak{L}^k$, the discrete de Rham sequence can be shown to be exact, thus equations (4.4) and (4.5) have stable approximations.

The discrete spaces for $H^* \mathfrak{L}^k$ must satisfy similar properties. Equations (5.2-5.3) hold in the dual sense, i.e. we write $\mathcal{P}_r^* \Lambda^k(\mathcal{T}_i^d) = \mathcal{P}_r^* \Lambda^k(\mathcal{T}_i^d) = \star (\mathcal{P}_r \Lambda^{d-k}(\mathcal{T}_i^d))$, and $d^* \mathcal{P}_r^* \Lambda^k(\mathcal{T}_i^d) \subset \mathcal{P}_r^{-*} \Lambda^{k-1}(\mathcal{T}_i^d) \subset \mathcal{P}_{r-1}^* \Lambda^{k-1}(\mathcal{T}_i^d)$. Furthermore, the coderivative $\mathfrak{d}^*$ imposes the inverted condition $\Lambda_h^{k+(n-e)}(\mathcal{T}^e) \subseteq \text{Tr}_{\Omega_i^{n-1}}^* \Lambda_h^k(\mathcal{T}^d)$ on boundaries.

## 6. Implications in terms of classical calculus

We take a moment to untangle the notation from Sections 3-5 in order to extract insight into modeling and discretization for the original physical problem.

Our initial task is to express simplest form of the mixed-dimensional Hodge Laplacian in terms of conventional notation. Considering for the moment $k = n$, the function spaces $H^* \mathfrak{L}^n$ and $H \mathfrak{L}^{n-1}$ correspond to $H_1$ scalars and $H(div)$ vectors on each dimension $d \geq 1$. For $d = 0$, only the scalars are defined. Furthermore, the term $\mathfrak{d}\mathfrak{a} \in \mathfrak{L}^{n+1} = 0$, and thus we arrive from (4.6) to the simpler problem

$$\mathfrak{b} = \mathfrak{K}\mathfrak{d}^*\mathfrak{a} \qquad \text{and} \qquad \mathfrak{d}\mathfrak{b} = \mathfrak{f}$$

(6.1)

In this case, the exterior derivative is the negative divergence plus jumps for each domain, while the codifferential is the gradient parallel to each domain, and the difference from boundaries perpendicular. As such, we arrive exactly at the model equations of Section 2, with the second choice of modeling assumption (2.10).

Turning our attention to the finite element spaces, the lowest order spaces for discretizing (4.4-4.5) are obtained by choosing $r_i^d = 1$ and $m_i^d = -$, from which we obtain piecewise constants for $\mathfrak{a}$ on all domains, while we obtain for $\mathfrak{b}$ the Nedelec 1$^{st}$ kind (div) − Raviart-Thomas − continuous Lagrange elements for domains with dimensions $d = 3,2,1$, respectively − all of the lowest order [10] (this method will be referred to as "Mixed 1$^{st}$ kind" in the next section). Interestingly, if we choose Nedelec 2$^{nd}$ kind (div) elements of lowest order for $d = 3$, equations (5.4-5.5) tell us that we should increase the order in the lower-dimensional domains, obtaining dG elements of order $n - d$ for pressure, with BDM (2$^{nd}$ order) − continuous Lagrange (3$^{rd}$ order) for fluxes in domains with $d = 2,1$. This is a new method resulting from the analysis herein. We refer to this method as "Mixed 2$^{nd}$ kind".

The mixed finite element discretization has the advantage of a strong conservation principle, and may be hybridized to obtain a cheaper numerical scheme (see [10] for a direct approach in this context, but also [6, 5] for direct constructions in the finite volume setting). Alternatively, we consider discretizing the Euler-variation of the unconstrained minimization problem, equations (4.3). The natural finite element spaces are $\mathcal{P}_r^{\mathfrak{m},*}\mathfrak{L}^n$, with $r_i^d = 1$ and $\mathfrak{m}$ does come into play, corresponding to $1^{\text{st}}$-order continuous Lagrange elements in all dimensions. From an engineering perspective, this formulation has been described in [18], we refer to this method as "Primal" in the next section.

A complimentary formulation is obtained by setting $k = 0$. In this setting, we note from the definition of $\mathfrak{L}^0$ in equation (3.1) that no variables exist in the sub-manifolds, thus the formulation consists of the model equations from Section 2, but with $k^d = 0$ for all $d < n$. Physically, this corresponds to barriers in the domain (see e.g. [7, 19]).

For $k = 1$ and $k = 2$, we obtain models for Maxwell's equations with thin material inclusions. As in the scalar case above, the two cases lead to different physical interpretation of the submanifolds [11].

## 7. Computational example

In order to illustrate the concepts discussed in the preceding sections, we will return to the example of $k = n$ and fractured porous media as a computational example, using the three numerical methods proposed for this problem in the previous section.

The example consists of the unit square with two fractures crossing through the domain, intersecting at a right angle, as illustrated in figures 3. We impose unit permeability in the surroundings, set the normal and tangential permeability of the fractures to 100 and assume the apertures of both fractures as $\epsilon = 10^{-3}$. The boundary conditions are chosen as zero pressure at the bottom and no-flux conditions on the sides. Moreover, a boundary pressure of one is imposed on the fracture crossing the top boundary. All computations were performed with the use of FEniCS [20].

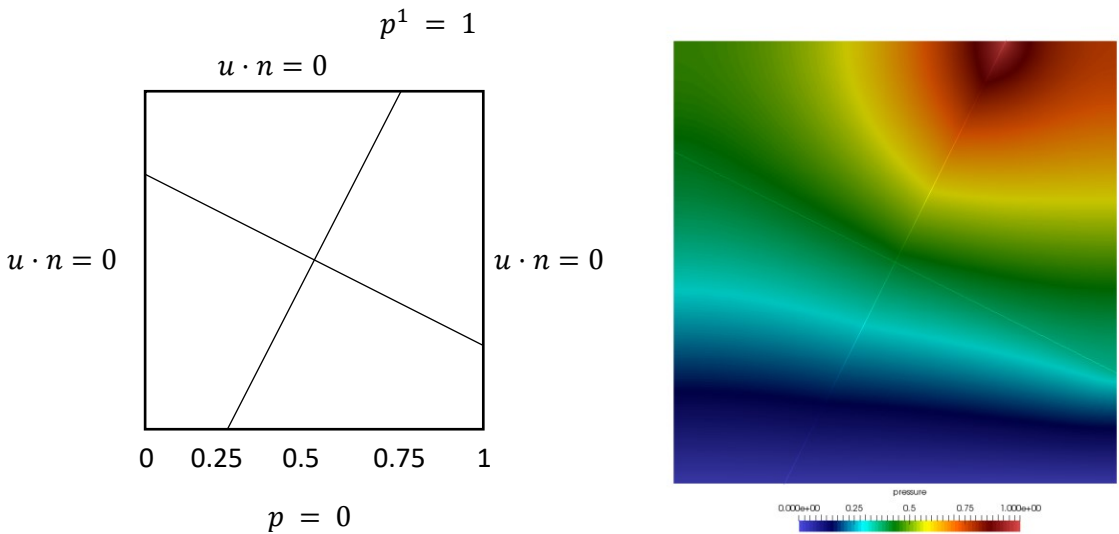

**Figure 3:** (Left) Domain of computation and associated boundary conditions. The pressure boundary condition is only imposed on the fracture pressure. (Right) Example of calculated solution (pressure).

The results show that all three methods are stable and convergent (Table 1). The relative errors and $L^2$-convergence rates after four consecutive refinements (identified by the characteristic grid size $h$) are given in the following table. Here, we compare the results to a fine-scale solution, obtained after a fifth refinement.

| Domain | Grid size | Primal Pressure | | Mixed 1st kind Pressure | | Flux | | Mixed 2nd kind Pressure | | Flux | |
|---|---|---|---|---|---|---|---|---|---|---|---|
| | $h$ | Error | Rate | Error | Rate | Error | Rate | Error | Rate | Error | Rate |
| $\Omega^0$ | $2^{-2}$ | 2.07e-02 | | 1.73e-02 | | N/A | N/A | 2.63e-03 | | N/A | N/A |
| | $2^{-3}$ | 7.67e-03 | 1.43 | 6.31e-03 | 1.46 | | | 8.81e-04 | 1.58 | | |
| | $2^{-4}$ | 2.66e-03 | 1.53 | 2.21e-03 | 1.52 | | | 2.89e-04 | 1.61 | | |
| | $2^{-5}$ | 8.45e-04 | 1.65 | 7.18e-04 | 1.62 | | | 8.99e-05 | 1.69 | | |
| | $2^{-6}$ | 2.15e-04 | 1.97 | 1.87e-04 | 1.94 | | | 2.26e-05 | 1.99 | | |
| $\Omega^1$ | $2^{-2}$ | 1.88e-02 | | 7.68e-02 | | 3.25e-02 | | 3.19e-03 | | 1.09e-02 | |
| | $2^{-3}$ | 7.00e-03 | 1.43 | 3.82e-02 | 1.01 | 1.48e-02 | 1.14 | 9.85e-04 | 1.70 | 4.47e-03 | 1.28 |
| | $2^{-4}$ | 2.54e-03 | 1.46 | 1.89e-02 | 1.01 | 6.32e-03 | 1.22 | 3.01e-04 | 1.71 | 1.84e-03 | 1.28 |
| | $2^{-5}$ | 9.57e-04 | 1.41 | 9.22e-03 | 1.04 | 2.49e-03 | 1.34 | 8.99e-05 | 1.74 | 7.44e-04 | 1.30 |
| | $2^{-6}$ | 3.23e-04 | 1.57 | 4.12e-03 | 1.16 | 7.82e-04 | 1.67 | 2.37e-05 | 1.92 | 2.61e-04 | 1.51 |
| $\Omega^2$ | $2^{-2}$ | 3.36e-02 | | 7.82e-02 | | 2.24e-01 | | 7.47e-02 | | 8.97e-02 | |
| | $2^{-3}$ | 1.23e-02 | 1.45 | 3.85e-02 | 1.02 | 1.37e-01 | 0.71 | 3.75e-02 | 1.00 | 5.31e-02 | 0.76 |
| | $2^{-4}$ | 4.25e-03 | 1.53 | 1.89e-02 | 1.02 | 8.21e-02 | 0.74 | 1.86e-02 | 1.01 | 3.16e-02 | 0.75 |
| | $2^{-5}$ | 1.36e-03 | 1.64 | 9.17e-03 | 1.05 | 4.75e-02 | 0.79 | 9.11e-03 | 1.03 | 1.87e-02 | 0.75 |
| | $2^{-6}$ | 3.60e-04 | 1.92 | 4.08e-03 | 1.17 | 2.47e-02 | 0.94 | 4.07e-03 | 1.16 | 1.04e-02 | 0.86 |

**Table 1:** Convergence rates for the three FE and MFEM discussed for the fracture problem in Section 6.

First, we observe that each method captures the intersection pressure well, with second order convergence over all. In the surroundings, the pressure convergence with second order for the primal formulation and first order for both mixed formulations, as expected. The Mixed 2nd kind method has higher-order elements in the fracture, and this is reflected in higher convergence rates for both pressure and flux.

## Acknowledgments

The authors wish to thank Gunnar Fløystad, Eirik Keilegavlen, Jon Eivind Vatne and Ivan Yotov for valuable comments and discussions on this topic. The research is funded in part by the Norwegian Research Council grants: 233736 and 250223.